\documentclass[11pt]{article}
\usepackage[noblocks]{authblk}
\usepackage{amsmath,amssymb}
\usepackage{subcaption}
\usepackage{inputenc}
\usepackage{graphicx}
\usepackage{float}
\usepackage{array,ragged2e}
\newcolumntype{M}[1]{>{\RaggedRight\arraybackslash}m{#1}}
\usepackage{amsthm}
\theoremstyle{plain}
\newtheorem{theorem}{Theorem}[section]

\theoremstyle{definition}

\newtheorem{example}[theorem]{Example}
\newtheorem{remark}[theorem]{Remark}

\usepackage{geometry}
\usepackage{float}
\usepackage[mathscr]{euscript}
\usepackage{cite}
\SetSymbolFont{operators}{bold}{OT1}{cmr}{bx}{n}
\SetSymbolFont{letters}{bold}{OML}{cmm}{b}{it}
\usepackage[font=scriptsize]{caption} 
\usepackage{etoolbox}

\makeatletter
\patchcmd{\maketitle}{\@fnsymbol}{\@alph}{}{}  
\makeatother
\title{\vspace{-3ex} \bfseries
	Dynamical behavior of a time-delayed infectious disease model with a non-linear incidence function under the effect of vaccination and treatment \vspace{-1ex}}
\author[a] {Sushil Pathak}
\author[b*]{G.Shirisha}
\author[a]{K.Venkata Ratnam}
\affil[a]{\small  Department of Mathematics, Birla Institute of Technology \& Science-Pilani,Hyderabad Campus, Jawahar Nagar, Hyderabad-500078, Telangana, India. p20210042@hyderabad.bits-pilani.ac.in, vrkota@hyderabad.bits-pilani.ac.in. }
\affil[b*]{\small Department of Mathematics, Stanley College of Engineering and Technology for Women, Abids, Hyderabad-500001, Telangana, India. deepasiri82@gmail.com}
\date{}
\begin{document}
	\maketitle \vspace{-3mm}
	\begin{abstract}
When an infectious disease propagates throughout society, the incidence function may rise at first due to an increase in pathogenicity and then decrease due to inhibitory effects until it reaches saturation. Effective vaccination and treatment are very helpful for controlling the effects of such infectious diseases. To analyze the impacts of these diseases, we proposed a new compartmental model with a generalized non-linear incidence function, vaccination function, and treatment function, along with time delays in the respective functions, which show how its monotonic features influence the stability of the model. Fundamental properties of a model, such as positivity, boundedness, and the existence of equilibria, are examined in this work. The basic reproduction number has been computed and correlative studies for local stability in view of the basic reproduction number have been examined at the disease-free and endemic equilibrium points. A delay-independent global stability result has been established, and to be more precise, we explicitly derived the result on global stability by restricting delay parameters within a very specific range. Furthermore, numerical simulations and some examples based on COVID-19 real-time data are pointed out to emphasize the significance of how the disease's dynamical behavior is characterized by various functions for controlling the spread of disease in a population and to justify the mathematical conclusions.
	\end{abstract}
	\textbf{Keywords}: SIR model, Equilibria, Boundedness, Uniqueness, Time Delays, Basic Reproduction Number, Local Stability, Lyapunov Functions, Global Stability, Simulations, and Numerical Results.\\
 \textbf{MSC}: 37D35, 34D20, 34D23, 93D05, 92D30, 92-10.
	\section{Introduction}
	Disease outbreaks and dissemination have been questioned and investigated for many years all around the world. The phrase communicable disease refers to conditions in which pathogens transfer from one person to another through a carrier or direct physical contact. Some infectious diseases like swine flu, chikungunya, dengue, plague, and malaria are transmitted by a medium such as water, air, or carriers like insects, flies, and mosquitoes \cite{BERETTA201887, LIU20153194, Mapder2019}. In recent times, one of the major global problems has been the infectious disease COVID-19, which emerged due to a virus named SARS-CoV-2. It has several mutants, such as Beta, Delta, Omicron, etc. Initially, this virus spread very rapidly among the population of Wuhan, China; after that, it spread worldwide and caused local epidemics in all countries \cite{BEKIROS2020109828, COOPER2020110298}. The main symptoms of the virus are fever, coughing, sneezing, sore throat, loss of taste and smell, etc., and this disease mainly affects people who are suffering from diseases like diabetes, stomach problems, lung problems, or heart snatching issues \cite{LIN2020211, NDAIROU2021110652}. Still, this virus has a significant effect on the human population. Recently, it has been noticed that the impact of the virus has been carried forward, and it has led to heart attacks, brain hemorrhages, and paralysis among the people of society. The aforementioned infectious diseases have caused millions of deaths worldwide. In these situations, mathematical models are a very important tool for describing the complex processes that allow infectious diseases to spread and to estimate disease transmission, recovery, fatalities, and other important parameters related to disease transmission in the population. Modeling techniques are also useful for understanding and predicting the probability and impact of a disease outbreak, as well as for providing crucial information to public health policymakers to minimize the effects of infection and how effectively to treat the disease present in society. \\[0.5mm]Researchers have studied many infectious disease models that are present in the literature \cite{dantas2018calibration,  li1995global, rao2019predictive, rao2013dynamic, ZHU2019115}. To analyze the impact of smallpox, Bernoulli suggested a mathematical modeling approach in 1760 \cite{Bernoulli17601}. In addition to this, H. Hamer outlined a mathematical model to estimate the periodicity of the measles epidemic in 1906 \cite{hamer1906epidemic}. The most straightforward model, which may be described in terms of disease state, is the SI epidemic model, where S stands for susceptible and I for infected. People in the SI model remain infected until they die. Additionally, the SI model is sometimes known as the SIS model, in which a person who has been infected and then healed becomes susceptible once more. In late 1927, W.O. Kermack proposed the first classical epidemic SIR model \cite{Kermack1927700}. The SIR model comprises three compartments susceptible (S), infected (I), and recovered (R). It is used to describe diseases in which immunity lasts for only a short period before vanishing \cite{kabir2019analysis}. In epidemiology, researchers have widely used SIR models, which help policymakers and healthcare professionals estimate the population's susceptibility, infection, and recovery rates to determine how much medical support is required for a particular disease. Further, they divided the population into numerous compartments and established several epidemic models depending on the factors necessary to examine the traits associated with the disease.
The SEIR, MSIR (maternally derived immunity), SEIS, and MSEIR \cite{li1995global, dantas2018calibration, de2021dynamics, giordano2020modelling, HATTAF2020109916, rao2013dynamic}, among other epidemic models, 
that have been investigated in the literature. \\[0.5mm] One of the important aspects in constructing the infectious disease model is the term disease incidence rate, which is described as the number of susceptible people who contract the disease per unit of time. Many infectious disease models operate under the presumption that the incidence rate is exactly proportional to the densities of the susceptible host along with the infection-carrier population. However, several biological parameters could cause non-linearity in the rate of disease spread. The incidence rate can also be modeled using more general functions. It is very crucial to take into account whether the functional form of the incidence rate can alter the dynamics of an epidemic. Therefore, some authors have proposed that the bilinear or nonlinear incidence rate be modified into various sorts of functions \cite{DONOFRIO2022112072, naresh2009stability}.\\[0.5mm]
 Now, when we used to study contagious disease transmission in epidemiology, coupled differential equations are frequently employed to anticipate the transmission of pathogenic diseases and explore the possible consequences of preventative interventions. More specifically, the study becomes more intriguing when there are unexpected delays in some particular functions of the model equations. As delay differential equations play a vital role in disease dynamics, the implementation of such time delays significantly enhances the complexity of the model \cite{perko2013differential}. There are several biological factors that may also lead to the introduction of different delays in epidemic models. The two most frequent causes of a delay in epidemiology are the infection's latency period in infected people and the infection's latency period in carriers. In these situations, it would take some time before the infection level in the infected host or the carriers approached a level sufficient to spread the virus in a population. Furthermore, in the early stages of a disease, treatment is an extremely effective technique for controlling the transmissible disease.
Treatment techniques may proceed in two ways: vaccination procedures or treatment through prescribed medicines. Hence, it is important to study the impact of vaccination and treatment through medicines to control the prevalence of the disease. Many existing models incorporate treatment rates hypothesized to be proportional to the number of infectious agents \cite{KABIR2019118, Zaman2017}. However, it is fairly obvious that communities have few medical resources available in the early stages of infectious disease for treating diseases. Since the delay may also happen to induce recovery in the infected host from the infectious disease due to boosting the immune system after treatment and vaccination in the human population. \\[0.5mm] 
As a result of the aforementioned considerations, we propose a SIR model, which was introduced in \cite{rao2015dynamic}, in which the terms were expressed as a coupled delay differential system that relates the linear, semi-linear, and non-linear terms together by forming a well-defined system of equations. In this model, authors have figured out the incidence rate along with the incidence function $f$ which relates susceptible individuals to infected individuals in the presence of incubation delay. In addition, they have also incorporated $v$ as a vaccination function and $p$ as a treatment function with a treatment delay. The dynamics of the model were given by, 
 \begin{align}
			x^{'} &= \hspace{0.2cm} a - b f(x,y) - dx - cv(x( t )) + \alpha z\hspace{0.1cm},\nonumber\\\vspace{0.2cm} 
			y^{'} &=\hspace{0.2cm} b_1f\Big( x(t - \tau) , \hspace{0.1cm}y\Big) - rp( y ) - d_1y\hspace{0.1cm},\nonumber\\\vspace{0.2cm} 
		z^{'}	&=\hspace{0.2cm} rp\Big( y(t - \delta)\Big) - \alpha z. \label{maineq}
		\end{align}
 Each parameter in this model has its own physical significance, which is going to be discussed in Section 2.
 So, in the context of the model, delays discussed in (\ref{maineq}), might be an unavoidable delay in vaccination function, which may arise due to a delay in the production of a suitable vaccine, a delay in the supply of vaccines, or, most importantly, the fear of side effects of vaccines in the human population. The influence of all three delays (i.e., delays in non-linear incidence function, treatment function, and vaccination function) together on disease prevalence is hardly studied in the literature. Henceforth, in this paper, we will examine the effect of vaccination, treatment, and incidence delays on the infected and recovered numbers. To verify the applicability of the model to the present scenario, we tried to analyze our model in real-time situations by applying three different real-time data sets. The first dataset is based on the omicron variant of COVID-19 in the state of Tamilnadu, India (December 2021–March 2022); the next two are the cumulative dataset (November 2020) of India COVID-19 and USA COVID-19 data \cite{dickson2022sqirv, wintachai2021stability}. \\ [1mm]The primary goal of this research article is as follows: First, construct model equations with respective delays in vaccination, infection, and treatment functions in Section 2. In Section 3, we will investigate the core aspects of the solutions, such as positivity, boundedness, and the existence-uniqueness of equilibria. The threshold value $R_0$ is determined for this particular model, and then local stability at endemic and disease-free equilibria is explored. By building appropriate Lyapunov functions, we analyze both delay-independent and delay-dependent conditions for global stability at the equilibria of the model in Section 4. In Section 5, numerical simulations of the system based on real-time data are implemented to validate theoretical results, followed by conclusions and future scope. \\ Based on the above discussion, we are directed toward describing our model.
	 \vspace{-4mm}
	\section{Description of the model} 
\begin{figure}[!ht]
    \centering
\includegraphics[height=6.4cm,width=11cm]{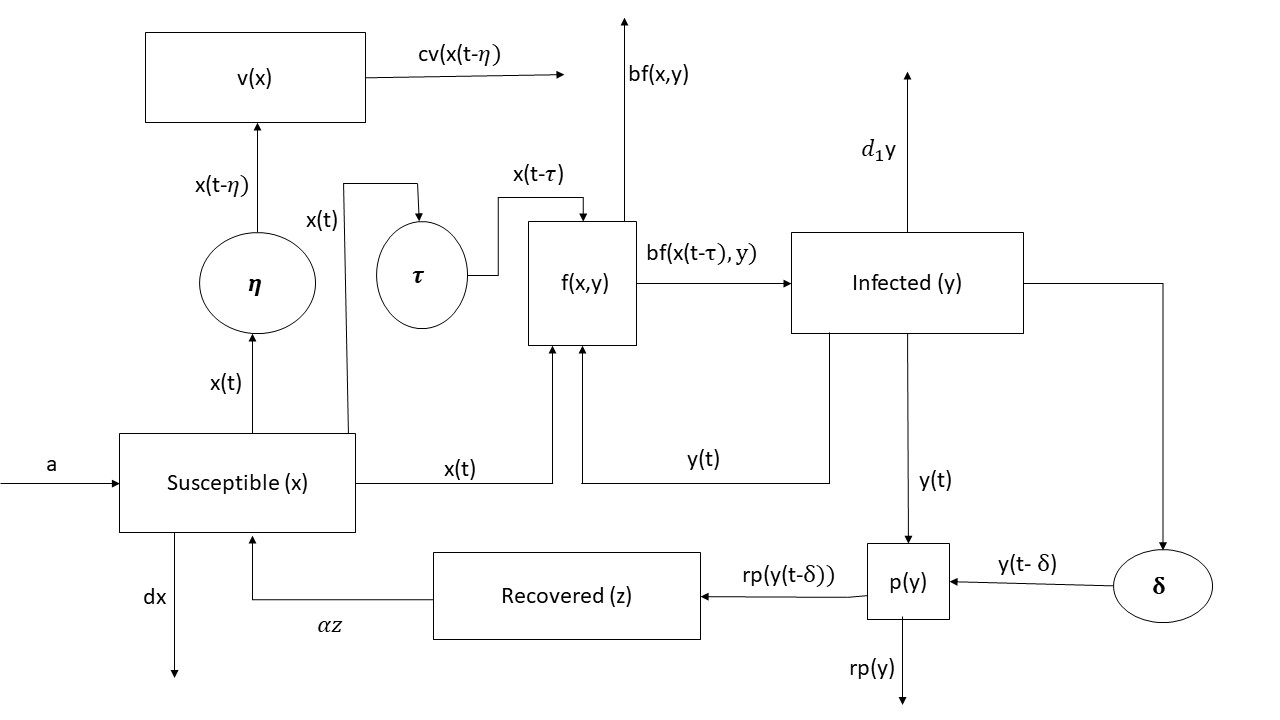}
    \caption{Block diagram of SIR model (\ref{eq1})}
    \label{block}
\end{figure}

\vspace{-0.5mm}
	The SIR model is an epidemiological model that predicts the potential prevalence of a contagious disease in a closed population over time. It is a compartmental model in which the model is classified into three distinct compartments named susceptible, infected, and recovered individuals. Because these models incorporate coupled equations that relate the number of susceptible people $x(t)$, the number of infected people $y(t)$, and the number of recovered people $z(t)$, they are named SIR models. Usually, these infectious diseases spread through carriers like mosquitoes, flies, infected people, etc. that are used to transmit the virus from infected people to susceptible people. As discussed in the above section, we are introducing a delay in the vaccination function in the SIR model (\ref{maineq}). Then, the following set of delayed differential equations provides the dynamics of the resulting model:\vspace{-2mm}
\begin{align}
x^{'} & = a - b f(x,y) - dx - cv\Big(x( t - \eta)\Big) + \alpha z,\nonumber\\
y^{'} & = b_1f\Big( x(t - \tau), y\Big) - rp( y ) - d_1y,\nonumber\\
z^{'} & = rp\Big( y(t - \delta)\Big) - \alpha z. \label{eq1} 
\end{align}\vspace{-2mm}
The parameters used in the model are listed in Table \ref{table1}.
\begin{table}[!ht]
    \small
    \centering
    \begin{tabular}{|c|c|}
		\hline
  $\textbf{Parameter}$ & \textbf{Biological significance} \\ 
  \hline
		$a$ & Evolution rate of the susceptible population  \\ 
		\hline
		$f$ & Infection function or non-linear incidence  \\  
		\hline
		$b$ & Contact rate of the infected with susceptible individuals in a population\\
		\hline
		$d$ & Natural death rate of the susceptible population\\
		\hline
		$v$ & Vaccination function\\
		\hline
		$c$ & Rate of effective vaccination\\
		\hline
		$b_1$ &  Effective rate of conversion of susceptible individuals into infected individuals \\
		\hline
		$p$ & Treatment function \\
		\hline
		$r$ & Rate of treatment \\
		\hline
		$\alpha$ & Rate at which rehabilitated individuals become susceptible again \\
		\hline
		$d_1$ & Death rate of infected individuals beyond the treatment\\
		\hline
  	$\eta$ & Delay in vaccination function\\
   \hline
   	$\tau$ & Delay in infection function\\
    \hline
    $\delta$ & Delay in treatment  function\\
    \hline
	\end{tabular}
 \caption{Description of model parameter} \label{table1}
 \end{table}\\\\
  Figure \ref{block} is a schematic representation of our model (\ref{eq1}). 
\\
 We are making the following assumptions:
\begin{itemize}
\item  $f(x,y) \geq 0$ for all $x$, $y$. Since $f(x,y)$ represents infection function, it should be positive for a well-defined biological system and its dynamical behavior. \item $f(0,y)=0$ for all $y$ because if there are nil susceptible cases in the population, then y does not make sense for infection function. \item $f(x,0)=0$ for all $x$ because if infected individuals are zero, then it will create a disease-free environment. \item $f(x,y)$ and $p(y)$ are always non-negative, continuously differentiable, and monotonically increasing for $x>0$, $y>0$. \item The incidence function, treatment function, and vaccination function can be chosen from the following classes of functions mentioned in the literature \cite{DONOFRIO2022112072, article1, naresh2009stability, rao2015dynamic, WANG20102390, ZHANG20081456}.\\[1mm]
$f(x,y)=cx^{p}y^{q}$ or $\dfrac{xy}{1+bx}$ or $\dfrac{xy}{1+ay},$\\[1.5mm]
$p(y)=k$ or $y$ or $\dfrac{y}{a + y}$ or $sinhy$ or $tanhy$,\\[1.5mm]
$v(x)=k$ or $x$ or $\dfrac{x}{a + x}$ or $sinhx$.
\item We also assume that the incidence function, vaccination function, and treatment function satisfy Lipschitz's condition. So, for positive constants $K_1, K_2, K_3, K_4$, $L_1, L_2, M_1, M_2$ we have,
\begin{align}
	 K_1 \vert x-\overline{x}\vert +K_2 \vert  y-\overline{y}\vert  &\leq \vert f(x,y)- f(\overline{x},\overline{y})\vert  \leq K_3 \vert  x-\overline{x} \vert +	K_4 \vert  y-\overline{y} \vert,  \nonumber\\
	 M_1 \vert  x-\overline{x} \vert &\leq \vert  v(x)-v(\overline{x}) \vert \leq M_2  \vert x-\overline{x}\vert,  \nonumber\\
	 L_1 \vert  y-\overline{y} \vert  &\leq \vert p(y)- p(\overline{y}) \vert  \leq L_2 \vert  y-\overline{y} \vert,   \label{lip1}
\end{align}
respectively for $x \neq \overline{x}$, $y \neq \overline{y}$. 	
Under the Lipschitz conditions (\ref{lip1}), by the hypothesis of delay differential equations, the solutions of the system (\ref{eq1}) exist.\end{itemize}
The model should be well-behaved to be applied to real-life circumstances. Therefore, the qualitative features of our model will be examined in the following section.
\section{Qualitative analysis of the model}
In the following section, we shall explore the fundamental characteristics of biological conditions. The solutions of the system behave adequately, and they should be non-negative because they express real quantities. Hence, these aspects are addressed in the next result:
		\subsection{Positivity and boundedness of solution}
  First, we will demonstrate the result that the solutions of the system (\ref{eq1}) are always non-negative.
		\begin{theorem}
		All solution(s) of the system (\ref{eq1}) with non-negative initial condition(s) is non-negative for all $t>0$ \cite{article2}.
  \end{theorem}
		\begin{proof}
		Since, we are assuming $f(x,y)\geq 0$, $f(0,y)=f(x,0)=v(0)=p(0)=0$ and initially there are no delays, i.e $\eta$ = $\tau$ = $\delta = 0$.\\
  From system (\ref{eq1}), we get; \\[1mm]
  At $x=0, x^{'} = a+ \alpha z \geq 0$, \\[1mm]
 At $y= 0, y^{'} = 0$, and \\[1mm]
At $z = 0, z^{'}	 = r p(y) \geq 0$ . 	
Thus, all the solutions of the system (\ref{eq1}) are non-negative for all $t>0$.
\end{proof}
 Generally, bounded solutions of the system seem to be predictable; hence, we will try to get a result on the boundedness of the solutions (\ref{eq1}).
\begin{theorem}
    The solution of the system (\ref{eq1}) is bounded provided the functions $f,v,p$ satisfy (\ref{lip1}).
\end{theorem}
\begin{proof} Consider $\Phi(t)$ as,
 \\ $\Phi(t) = \displaystyle{ \vert x(t) \vert + \vert y(t) \vert +\vert  z(t) \vert -c M_1 \int_{t-\eta}^{t} \vert x(s) \vert ds + b_1 K_3 \int_{t-\tau}^{t} \vert x(s) \vert ds+r L_2 \int_{t-\delta}^{t} \vert y(s) \vert ds}$, \\
 for all $t \geq 0$. Then, the dini derivative along the solution of (\ref{eq1}) is, 
    \begin{align*}
        D^+ \Phi(t)& \leq a- b\vert f(x,y) \vert-d \vert x \vert -c \vert v(x(t-\eta)) \vert + \alpha \vert z(t) \vert +b_1 \vert f(x(t-\tau),y) \vert-d_1 \vert y (t) \vert \\&- r \vert p(y(t)) \vert +  r \vert p(y(t-\delta)) \vert - \alpha \vert z(t) \vert-c M_1 \vert x(t) \vert+c M_1 \vert x(t-\eta) \vert+ b_1 K_3\vert x(t) \vert\\&- b_1 K_3\vert x(t-\tau) \vert+rL_2 \vert y(t) \vert- r L_2 \vert y(t-\delta) \vert,
    \end{align*}
    Using the conditions given in (\ref{lip1}), we will get;
         \begin{align*}
        D^+ \Phi(t) & \leq a-b K_1 \vert x(t) \vert -b K_2 \vert y(t) \vert -d \vert x(t) \vert -cM_1 \vert x(t) \vert + b_1 K_3\vert x(t) \vert +b_1 K_4 \vert y(t) \vert \\&-d_1 \vert y(t) \vert -r L_1 \vert y(t) \vert+ r L_2\vert y(t) \vert,\\
        & \leq a- (b K_1+d+cM_1- b_1 K_3) \vert x(t) \vert -(b K_2-b_1 K_4+d_1+r L_1-r L_2) \vert y(t) \vert,\\
        & \leq a-\gamma \Phi(t),
    \end{align*}
    where, $\gamma=\min \Big \{b K_1+d+cM_1- b_1 K_3, b K_2-b_1 K_4+d_1+r L_1-r L_2  \Big \}$.\\
    By employing the comparison theorem \cite{comparisonthm1}, we get $\Phi(t) \leq\max \Big\{\dfrac{a}{\gamma}, \Phi(0) \Big\}$. Thus, the solutions of the system (\ref{eq1}) are bounded.
\end{proof}		
	Thus, solutions of the model (\ref{eq1}) are positive and bounded. However, it is essential to understand how the solutions will perform over time. As the concept of the dynamical behavior of a biological system maintaining equilibrium relates to its attractors, which characterize its long-term behavior and it is crucial to examine the factors that contribute to the stability of equilibria. So, in this context, we shall discuss the results of the existence-uniqueness of equilibria and local as well as global stability in the next sections.
		\subsection{Equilibria}
		\vspace{-0.5mm}An equilibrium of a dynamical system is a constant value of the state variables that does not fluctuate over time \cite{perko2013differential}. In our model, there will be two types of equilibria: disease-free equilibria and endemic equilibria. Disease-free equilibria is defined as the point when there is no disease infection present in the population $(x^*, 0, 0)$, where $x^{*} > 0$. Endemic equilibria are defined as the point when there is a disease outbreak consistently present in the population, that is, $(x^*, y^*, z^*)$, where $x^{*} > 0$, $y^{*} > 0$, and $z^{*} > 0$. Throughout the rest of the study, $(x^*, y^*, z^*)$ signifies either of these equilibria unless otherwise specified.\\
 As discussed above, equilibrium stability is crucial to understanding how the model eventually behaves. We first deduce a result that gives the conditions under which a unique equilibrium point exists for our system (\ref{eq1}).
		\begin{theorem}
	Let the infection function $f(x, y)$, vaccination function $v(x)$, and treatment function $p(y)$ satisfy the Lipschitz conditions (\ref{lip1}). Then the system has a unique equilibrium solution if the parameters adhere to the following conditions:
 \begin{align}
     A &=  \Bigg(\dfrac{b_1 K_3}{d_1}-\dfrac{b K_3}{d} -\dfrac{c M_1}{d}\Bigg)< 1,\nonumber\\
		   B &= \Bigg(\dfrac{b_1 K_4}{d_1}-\dfrac{b K_4}{d}+\dfrac{rL_2}{\alpha}-\dfrac{rL_2}{d_1} \Bigg)< 1,\nonumber\\
		  C &=  \dfrac{\alpha}{d} <1. \label{eqcond}
 \end{align}
		\end{theorem}
		\begin{proof}
		Let us consider $\textbf{D} = \{X \in R^{3} \lvert \lvert X \rvert< K\}\subseteq R^{3}\}$, where $\textbf{D}$ is invariant in the context of the delay differential equation.\\[2mm]
		Define $\lvert X \rvert = \lvert x \rvert + \lvert y \rvert + \lvert z \rvert$, \\[2mm]
		We can write the system as $X^{'} = PX $, where $X = (x, y, z)$,\\[1mm] $P(X) = P(x, y, z)$ =  $\Big(F(x, y, z), G(x, y, z), H(x, y, z)\Big)$,\\[1.5mm]
		$F(x,y,z)$ = $\dfrac{a - bf(x,y) - cv(x(t - \eta)) + \alpha z}{d}$,\\[1mm] 
			$G(x,y,z)$ = $\dfrac{b_1f(x(t - \tau ),y) - rp(y)}{d_1}$,\\[1mm]
	$H(x,y,z)$ = $\dfrac{r p(y(t - \delta))}{\alpha}$.\\[1.5mm]
 Now, let us assume $X,\Bar{X} \in \textbf{D}$,
		\begin{align*}
		\lvert P(X) - P(\Bar{X})\rvert 		   & =\lvert(F(x, y, z) - F(\Bar{x},\Bar{y},\Bar{z})\rvert + \lvert (G(x, y, z) - G(\Bar{x},\Bar{y}, \Bar{z})\rvert \\&+ \lvert(H(x, y, z) - H(\Bar{x}, \Bar{y}, \Bar{z})\rvert,\\[2mm]
		    \leq & \Bigg(\dfrac{b_1K_3}{d_1} -\dfrac{bK_3}{d} -\dfrac{cM_1}{d}\Bigg)\lvert x - \Bar{x}\rvert + \Bigg(\dfrac{b_1K_4}{d_1} - \dfrac{bK_4}{d}  +\dfrac{rL_2}{\alpha} -\dfrac{rL_2}{d_1} \Bigg)\lvert y  - \Bar{y}\rvert \\&+ \dfrac{\alpha}{d}\lvert z - \Bar{z}\rvert.
		\end{align*}
		   By taking $\Bar{A}= Max(A, B, C) $, we get;\\[2mm]
		 $\lvert P(X) - P(\Bar{X})\rvert \leq \Bar{A}\Big[\lvert x - \Bar{x} \rvert + \lvert y - \Bar{y}\rvert + \lvert z - \Bar{z}\rvert\Big]
		 \leq \Bar{A} \lvert X - \Bar{X} \rvert.\\[1mm]$
   Hence, by the contraction mapping theorem, there exists a unique point $(x^*, y^*, z^*)$ such that $P(x^*, y^*, z^*)$ = $(x^*, y^*, z^*)$. Thus, the system (\ref{eq1}) has a unique solution.
\end{proof}
Now, our assumption for the conditions (\ref{eqcond}) holds, i.e., there exists a unique equilibrium point $(x^*,y^*,z^*)$ for the system (\ref{eq1}).\\[1mm]
Next, we are going to derive the basic reproduction number for the model (\ref{eq1}), which is an essential aspect of epidemic models. It is frequently used to assess the severity of an epidemic outbreak since it measures a disease's propensity for transmission. Using a basic reproduction number, we will establish the model's local stability at disease-free equilibrium. We also try to establish the local stability of the model at endemic equilibrium by restricting the parameters.
	\subsection{Basic Reproduction number (BRN)}
In epidemiology, the basic reproduction number is a measurement that is used to calculate the number of secondary infections induced by an infectious person moving in a susceptible population. It is denoted by $\textbf R_0$. In a disease model, the importance of $\textbf R_0$ is depicted by its values, and it also represents the strength of disease transmission within a community. If $\textbf R_0$ $>$ 1 then each person will transmit the infection to more than one person, and the disease transmission rate will increase. On the other hand, if $\textbf R_0$ $<$ 1, then each infected individual will transmit the infection to less than one person. As a result, the disease will disperse throughout the population.\\[1mm]
  We will employ the next-generation matrix methodology for the calculation of the basic reproduction number. In the next-generation matrix, there will be terms like F and V, which are described in \cite{article3}. It is usually used as a threshold parameter for contagious disease and specified as $R_0$ = $\rho(FV^{-1})$, where $\rho$ is the spectral radius operator.\\[1mm]
In the model formulation, there is only one infected compartment, as shown in (\ref{eq1}). Now our job is to generate an F and V matrix where F is defined as the rates of flows from uninfected to infected classes and V is the rates of all other flows to and from the infected class at disease-free equilibrium \cite{Singh_2022}.\\To calculate the threshold $\textbf R_0$, we have to estimate the spectral radius of the next-generation matrix \cite{article3} $T = FV^{-1}.$\\[2mm]
Hence, for our model, value of $F$ and $V$ are as follows:\\[2mm]
	  $F = \begin{bmatrix}
	      \dfrac{\partial f_1}{\partial y} \end{bmatrix} = b_1f_y( x^{*}, 0) $\hspace{2cm}  $ V = \begin{bmatrix}
	     \dfrac{\partial v_1}{\partial y}\end{bmatrix} =  rp^{'}(0) + d_1;$ \\[3mm]
      where, $ f_1 = \begin{bmatrix}
    b_1 f_y (x(t - \tau),y) 
 \end{bmatrix}$ and $v_1 = \begin{bmatrix}
    rp(y) + d_1y
 \end{bmatrix}. $\\[2mm]
 The next-generation matrix is  $T = \begin{bmatrix}
     \dfrac{b_1 f_y (x^{*},0)}{rp^{'}(0) + d_1} 
 \end{bmatrix}$.\\[2mm]
 Now, to obtain the reproduction number $\textbf R_0$, we have to figure out the largest eigenvalue of $T$.
 \begin{center}
$det(\lambda I - T) = \begin{bmatrix}
    \lambda -  \dfrac{b_1 f_y (x^{*},0)}{rp^{'}(0) + d_1} 
 \end{bmatrix}$, 
 \end{center} Then,
	 \begin{center}
 $R_0 = \vert \rho(FV^{-1})\vert = \dfrac{b_1 f_y (x^{*},0)}{ rp^{'}(0) + d_1 }.$\end{center}
	 
Next, we try to establish the local stability of the model (\ref{eq1}) at endemic and disease-free equilibria.
	\subsubsection{Local stability}
 The local stability of the system reveals the stability of the system along with the behavior of the solution around equilibrium points. Now, we shall derive the criteria required for the local stability of the model at disease-free and endemic equilibrium points.
  \begin{theorem}
      For \textbf{$R_0$} $>$ 1, the system (\ref{eq1}) will be unstable at disease-free equilibrium. For \textbf{$R_0$} $<$ 1, it will be locally stable provided the condition $\dfrac{r}{\alpha}e^{-\lambda \delta}  \dfrac{dp(0)}{d y_\delta}< 1 $ holds true. \label{local2}
  \end{theorem}
  \begin{proof}
  Characteristic equation of (\ref{eq1}) at $(x^*,0,0)$ is, \begin{align*}
\begin{vmatrix}
 \lambda - T_1 & D_1 & \alpha\\
0 & \lambda - F_1 & 0\\
0 & 0 & \lambda - H_1
\end{vmatrix}\hspace{.3cm}=\hspace{.3cm} 0,
\end{align*}
where, $T_1 = - b f_x(x^*,0) - d -ce^{-\lambda\eta} \,\,v'(x^*)$, 
 $D_1 = -b f_y(x^*,0),$\\[2mm]
 $F_1 =   b_1 f_y(x^*,0) - rp^{'}(0) - d_1 + b_1e^{-\lambda \tau}  f_{x_\tau}(x^*,0), $
  $H_1 = - \alpha + r e^{-\lambda\delta}\,\,  \dfrac{dp (0)}{d y_{\delta}}$.\\[2mm]
  As we know, the system will be unstable even if one of the eigenvalues has a positive real part. \\[2mm] Let us take an eigenvalue $F > 0$, we get;\\[2mm]
   $\Bigg( b_1 f_y(x^*, 0) - rp^{'}(0) - d_1 + b_1 f_{x_\tau}(x^*, 0)\Bigg) > 0$\\[2mm]
   Since, $b_1 > 0,  f_{x_\tau}(x^*,0) > 0, \implies b_1 f_{x_\tau}(x^*,0) > 0$. \\[2mm]
   Therefore, $ b_1 f_y(x^*, 0) - r p^{'}(0) - d_1 > 0$, $\implies$
   $ b_1 f_y(x^*, 0) > rp^{'}(0) + d_1$,\\[2mm]
   $\dfrac{ b_1 f_y(x^*, 0)}{rp^{'}(0) + d_1}  > 1 \implies
   R_0 > 1$\\[3mm]
   Thus, the system is unstable if $R_0 > 1$.\\[2mm]
   For the system to be stable, all the eigenvalues must be negative. \\[2mm]Clearly, $T_1<0$.\\[2mm]
   By taking $F_1<0$, we get $b_1 f_y(x^*,0) - rp^{'}(0) - d_1 + b_1e^{-\lambda \tau}  f_{x_\tau}(x^*,0)<0$.\\[2mm]
   Since, $b_1e^{-\lambda \tau}  f_{x_\tau}(x^*,0)>0$ 
 $\implies b_1 f_y(x^*,0) - r p^{'}(0) - d_1 <- b_1e^{-\lambda \tau}  f_{x_\tau}(x^*,0)<0 $.\\[2mm]
   Therefore, $b_1 f_y(x^*,0) - rp^{'}(0) - d_1 <0$ $\implies R_0<1$.\\[2mm]  
   By $H_1<0$, we get $re^{-\lambda \delta} \dfrac{dp(0) }{d y_\delta}< \alpha $.\\
   Thus, the system is stable for $R_0 < 1$ and $\dfrac{r}{\alpha}e^{-\lambda \delta}  \dfrac{dp(0)}{d y_\delta}< 1 $.  
	\end{proof}
 Next, we try to deduce the conditions for local stability at endemic equilibrium.
 \begin{theorem}
 Endemic equilibrium point $(x^*, y^*, z^*)$ of the system (\ref{eq1}) is locally asymptotically stable if
 \begin{align}
     \dfrac{b_1}{rp'(y^*) + d_1}\Bigg( f_y(x^*, y^*)+e^{-\lambda \tau} f_{x_\tau}(x^*, y^*)\Bigg)<1 \,\,\,\, \mbox{and} \,\,\, \dfrac{r}{\alpha} e^{-\lambda\delta}  \dfrac{dp(y^*)}{d y_\delta}< 1. \label{loc1}
 \end{align}  \label{local1}
 \end{theorem}
\begin{proof}
Characteristics equation of system (\ref{eq1}) is given below:
\begin{align*}
\begin{vmatrix}
 \lambda - T_2 & D_2 & \alpha\\
0 & \lambda - F_2 & 0\\
0 & 0 & \lambda - H_2
\end{vmatrix} \hspace{.3cm}=\hspace{.3cm} 0
\end{align*}
where,\vspace{0.3cm}\\
    $T_2 = - bf_ x(x^{*},y^{*}) - d -ce^{-\lambda\eta} \,\,v'(x^{*})$, $D_2 = -b f_ y(x^{*},y^{*})$,\\[2mm] $F_2 =   b_1 f_{y} (x^{*}, y^{*})- rp^{'}(y^{*}) - d_1 + b_1e^{-\lambda\tau} f_{x_\tau}(x^{*}, y^{*})$,  $H_2 = - \alpha + r e^{-\lambda\delta}\,\,  p' (y^{*}) \vspace{0.3cm}$.\\   Eigenvalues are given by; $\lambda_1 = T_2$, $\lambda_2 = F_2$, $\lambda_3 = H_2$.\\[1.5mm]   Next, we know that if all the eigenvalues have negative real parts, then the system will be locally stable \cite{rao2015dynamic}.\\[2mm]  Clearly, $T_2<0$; since, $b f_x(x^*,y^*) + d +ce^{-\lambda\eta} \,\,v'(x^*)>0$.\\[2mm]  By taking $F_2<0$ and $H_2<0$, we get; $\dfrac{b_1}{rp'(y^*) + d_1}\Bigg(f_y(x^*, y^*)+e^{-\lambda \tau}  f_{x_\tau}(x^*, y^*)\Bigg)<1$ and $ \dfrac{r}{\alpha} e^{-\lambda\delta}  \dfrac{dp(y^*)}{d y_\delta}< 1 $, \\[2mm]which are our assumptions (\ref{loc1}). This completes the proof.
  \end{proof}
	\begin{remark} 
The key benefit of Theorem \ref{local2} is its ability to determine the model's instability at disease-free equilibrium solely based on the values $R_0$. On the contrary,  we need to restrict the parameters along with the values of $R_0$ for the model to acquire stability. Besides, Theorem \ref{local1} provides parametric constraint to validate local stability criteria at the endemic equilibrium point.
Since our model is more generalized than the model explored in \cite{rao2015dynamic}, the basic reproduction number derived here and Theorems \ref{local2} and \ref{local1} (results on local stability) can also be applied to the model \cite{rao2015dynamic}. 
	\end{remark}
One of the key concerns in epidemiology is the exploration of the global stability of endemic equilibrium in a mathematical model of contagious disease \cite{martorano2015equilibrium}. In addition to local stability, it is necessary to analyze the global stability of the model as it determines how the model ultimately behaves regardless of the initial circumstances, whether the infection eventually disappears, remains under control, or remains at high critical levels. So, in the following section, we will establish the results on the global stability of the model in two ways: delay-independent stability and delay-dependent stability.
	 \section{Global stability}
The global aspects of systems in the sense of the model's stability have been quite complicated in disease transmission models. Researchers have proposed a number of techniques to deal with the stability of equilibrium points, including the use of geometrical approaches, monotone dynamical systems, along with the Lyapunov functions method \cite{article5,liao2007stability, TIAN201131}. The global stability of the endemic equilibrium is particularly interesting whenever describing the behaviors of epidemic models on complex networks \cite{naresh2009stability}. Here, our approach is to show the global stability of the model by generating appropriate Lyapunov functions and extrapolating some parametric conditions. 

   Basically, an equilibrium point is a stationary solution of the system; hence, from (\ref{eq1}), we get;
\begin{align}
	a-b f(x^{*},y^{*})-d x^{*}-cv(x^{*})+\alpha z^{*}&=0, \nonumber\\
	b_1 f(x^{*},y^{*})-r p(y^{*})-d_1 y^{*}&=0,\nonumber\\
	r p(y^{*})-\alpha z^{*}&=0. \label{equil}
\end{align}
 From equations (\ref{eq1}) and (\ref{equil}), we get;
\begin{align}
	(x-x^{*})^{'} &=-b (f(x,y)-f(x^{*},y^{*}))-d (x- x^{*})-c(v(x(t-\eta))-v(x^{*}))+\alpha (z- z^{*}), \nonumber\\
	(y-y^{*})^{'} &=b_1 (f(x(t-\tau),y)- f(x^{*},y^{*}))-r (p(y) - p(y^{*}))-d_1(y- y^{*}),\nonumber\\
	(z-z^{*})^{'} &=r(p(y(t-\delta))- p(y^{*}))-\alpha(z - z^{*}).\label{equil1}
\end{align}
Now, let us discuss the equilibrium solution for global stability by restricting the parameters (delay-independent) and the delays of the model (delay-dependent). \\\subsection{Delay independent stability}Here, we restrict the parameters of the model under which the solutions of the system reach equilibria globally.
\begin{theorem}
	The positive equilibrium (endemic equilibrium) solution $(x^*, y^*, z^*)$ is globally asymptotically stable (in the absence of time lags), provided the functions $f$, $v$, and $p$ satisfy the Lipschitz conditions and the parameters of the system satisfy the conditions $b K_1 -b_1 K_3 +d+c M_ 1>0$ and $b K_2-b_1 K_4+ r  L_1+r  L_2+d_1 >0$. \label{indep}
\end{theorem}
\begin{proof}
Let the Lyapunov function be $U_1=\vert x-x^{*} \vert+ \vert y-y^{*}\vert+\vert z-z^{*}\vert$.\\ Now, the dini derivative along the solutions, using (\ref{equil1}) and (\ref{lip1}) is,
\begin{align}
	& D^{+}U_1\leq -b K_1 \vert x-x^{*}\vert- b K_2 \vert  y-y^{*} \vert-d \vert x- x^{*} \vert-c M_1 \vert x(t-\eta)-x^{*}\vert+b_1 K_3 \vert  x(t-\tau)-x^{*} \vert\nonumber\\
	& +b_1 K_4 \vert  y-y^{*} \vert-r  L_1 \vert  y-y^{*} \vert -d_1 \vert y- y^{*} \vert + rL_2 \vert  y(t-\delta)-y^{*} \vert. 
\end{align}	
Let $U_2=b_1 K_3 \displaystyle\int_{t-\tau}^{t} \vert  x(u)-x^{*} \vert du - c M_1 \int_{t-\eta}^{t} \vert  x(u)-x^{*} \vert du + rL_2 \int_{t-\delta}^{t} \vert  y(u)-y^{*} \vert du$.\\
Then,
\begin{align}
	D^{+}U_2&\leq b_1 K_3 \vert  x(t)-x^{*} \vert-b_1 K_3 \vert  x(t-\tau)-x^{*} \vert-c M_1 \vert x(t)-x^{*}\vert+c M_1 \vert x(t-\eta)-x^{*}\vert \nonumber\\&+ rL_2 \vert  y(t)-y^{*} \vert - rL_2 \vert  y(t-\delta)-y^{*} \vert. 
\end{align}
Let $U=U_1+U_2$,\\[1mm] As a result, 
\begin{align}
	 D^{+}U &\leq - \Bigg(\Big(b K_1 -b_1 K_3 +d+c M_1\Big)\vert x-x^{*}\vert+\Big(b K_2-b_1 K_4+ r  L_1+r  L_2+d_1  \Big)\vert  y -y^{*} \vert \Bigg)<0.
\end{align}
By the standard argument, we can say that $x(t)\longrightarrow x^*$,  $y\longrightarrow y^*$ and $z\longrightarrow z^*$ as $t \longrightarrow \infty.$ Therefore, the equilibrium $(x^*,y^*, z^*)$ is globally asymptotically stable.
\end{proof}

Next, we will construct a result based on delay-dependent stability, wherein we restrict the delay parameters for the global stability of the model.
\subsection{Delay dependent stability}
Here, we impose the parametric conditions on the system to get the delay range in view of global stability. \\For simplicity, we assume,\\[1mm]
$X = x -  x^*$, $Y = y -  y^*$, $Z = z -  z^*$, \\[1mm]$F(X,Y) = f(x, y) - f(x^{*},y^{*})$, $P(Y) = p(y) - p(y^*)$, $V(X) = v(x) - v(x^*)$.\\[1mm]
Then, the system (\ref{eq1}) will be transformed to,
\begin{align}
			X^{'} & = a - b F(X,Y) - dX - cV\Big(X( t - \eta)\Big) + \alpha Z, \nonumber\\ 
			Y^{'} &= b_1F\Big( X(t - \tau) , \hspace{0.1cm}Y\Big) - rP( Y ) - d_1Y,\nonumber\\
		Z^{'}&= rP\Big( Y(t - \delta)\Big) - \alpha Z. \label{eq11}
		\end{align}
The Lipschitz conditions (\ref{lip1}) can be written as,
\begin{align}
	 K_1 \vert X\vert +K_2 \vert  Y \vert  &\leq \vert F(X,Y)\vert  \leq K_3 \vert X  \vert +	K_4 \vert  Y \vert  \nonumber,\\
	 M_1 \vert  X \vert &\leq \vert  V(X) \vert \leq M_2  \vert X\vert  \nonumber,\\
	 L_1 \vert  Y\vert  &\leq \vert P(Y) \vert  \leq L_2 \vert  Y \vert.  \label{Lips2}
\end{align}
\begin{theorem}
	Suppose the functions of the system (\ref{eq11}) satisfy (\ref{Lips2}) and the parameters satisfy the conditions, 
	\begin{align*}
		A_1&= \Big(bK_1 - b_1K_3 + d_1 + cM_1 \Big)>0,
		A_2= \Big(bK_2 - b_1K_4 + rL_1 - rL_2 +d_1 \Big)>0,\\
            B_1&= \Big( cbM_1K_1 + cdM_1 + c^{2}M_1M_2\Big),
            B_2= \Big( b_1cM_1K_2\Big),
            B_3= \Big( cM_1 \alpha \Big)>0,\\
            C_1&= \Big(b_1K_3^{2}b + b_1dK_3 +b_1K_3cM_2 \Big),
            C_2= \Big(b_1bK_3K_4 \Big),
            C_3= \Big(b_1 K_3 \alpha \Big),\\
            D_1&= \Big( rL_2b_1K_3\Big),
            D_2= \Big(rL_2b_1K_4 + r^{2}L_2^{2} + rL_2d_1 \Big),
	\end{align*}
	then, the equilibrium point $(x^{*}, y^{*}, z^{*})$ of system is globally asymptotically stable for $0<\eta<r$, $0<\tau<s$ and $0<\delta<q$, where, $r=Min\Big\{\dfrac{A_1}{B_1}, \dfrac{A_2}{B_2}\Big\}$, $s=Min\Big\{\dfrac{A_1}{C_1}, \dfrac{A_2}{C_2}, \dfrac{A_1 B_3}{C_3 B_1}, \dfrac{A_2 B_3}{C_3 B_2}\Big\}$ and $q=Min\Big\{\dfrac{A_1}{D_1}, \dfrac{A_2}{D_2}\Big\}.$
    \label{dep}
\end{theorem} \vspace{-2mm}
\begin{proof}
Let the Lyapunov function be $E_1 = \lvert X \rvert + \lvert Y \rvert + \lvert Z \rvert$, \vspace{-2mm}
	\begin{align}
 \mbox{As we know,}\,\,\,\,
	X(t - \eta)& = X(t) - \displaystyle\int_{t-\eta}^{t} X^{'}(S)dS,\,\,\,\,
	X(t - \tau)& = X(t) - \displaystyle\int_{t-\tau}^{t} X^{'}(S)dS,\nonumber\\
	Y(t - \delta)& = Y(t) - \displaystyle\int_{t-\delta}^{t} Y^{'}(S)dS. \label{ddc1}
	\end{align}
Furthermore, by taking the dini derivative along with the solutions of (\ref{eq11}) by using (\ref{Lips2}) and (\ref{ddc1}) is,
 \begin{align*}
&D^+ E_1 \leq -(b K_1 + d + c M_1 - b_1K_3)\lvert X(t) \rvert + (-b K_2 + b_1K_4 - r L_1 - d_1 + rL_2) \lvert Y(t) \rvert \\& - c M_1 b K_1 \displaystyle\int_{t-\eta}^{t} \lvert X(S) \rvert dS - c M_1 b K_2 \displaystyle\int_{t-\eta}^{t} \lvert Y(S) \rvert dS- c M_1d \displaystyle\int_{t-\eta}^{t} \lvert X(S) \rvert dS  \\&- c^{2}M_1M_2 \displaystyle \int_{t-\eta}^{t}\lvert X(S-\eta) \rvert dS - cM_1\alpha \displaystyle \int_{t-\eta}^{t} \lvert Z(S) \rvert dS+ b_1bK_3^{2} \displaystyle \int_{t-\tau}^{t} \lvert X(S) \rvert dS \\&+ b_1 b K_3 K_4 \displaystyle\int_{t-\tau}^{t}\lvert Y(S) \rvert dS  + b_1d K_3 \displaystyle\int_{t-\tau}^{t}\lvert X(S) \rvert dS + b_1 b K_3 c M_2 \displaystyle\int_{t-\tau}^{t}\lvert X(S-\eta) \rvert dS+ \\&b_1 K_3\alpha \displaystyle\int_{t-\tau}^{t}\lvert Z(S) \rvert dS + b_1 r L_2 K_3 \displaystyle \int_{t-\delta}^{t} \lvert X(S-\tau) \rvert dS+  b_1 rL_2 K_4 \displaystyle\int_{t-\delta}^{t} \lvert Y(S) \rvert dS \\&+ r^{2} L_2^{2} \displaystyle\int_{t -\delta}^{t}\lvert Y(S) \rvert dS + d_1 r L_2\displaystyle\int_{t-\delta}^{t}\lvert Y(S) \rvert dS.
\end{align*}
Now, let us consider $E_2$ as,
\begin{align*}
E_2&= -c M_1 b K_1 \displaystyle \int_{t-\eta}^{t} dS \displaystyle \int_{s}^{t}\lvert X(U) \rvert dU -c M_1 b K_2 \displaystyle \int_{t-\eta}^{t} dS \displaystyle \int_{s}^{t} \lvert Y(U) \rvert dU -c M_1 d \displaystyle \int_{t-\eta}^{t} dS \displaystyle \int_{s}^{t} \lvert X(U) \rvert dU\\&-c M_1 \alpha \displaystyle \int_{t-\eta}^{t} dS \displaystyle \int_{s}^{t} \lvert Z(U) \rvert d U-c^{2} M_1 M_2 \Big( \displaystyle \int_{t-\eta}^{t}dS \displaystyle \int_{s}^{t}\lvert X(U-\eta) \rvert dU + \eta \displaystyle \int_{t-\eta}^{t} X(S) dS \Big) \\&+ b_1 bK_3^{2} \displaystyle \int_{t-\tau}^{t} dS \displaystyle \int_{s}^{t}\lvert X(U) \rvert + b_1 b K_3 K_4 \displaystyle \int_{t-\tau}^{t}dS \displaystyle \int_{s}^{t} \lvert Y(U) \rvert  + b_1 dK_3 \displaystyle \int_{t-\tau}^{t} dS \displaystyle \int_{s}^{t} \lvert X(U) \rvert \\&+ b_1 K_3 \alpha \displaystyle \int_{t-\tau}^{t} dS \displaystyle \int_{s}^{t} \lvert Z(U) \rvert +b_1 K_3 c M_2 \Big(\displaystyle \int_{t\tau}^{t} dS \displaystyle \int_{s}^{t} \lvert X(U-\eta) \rvert dU + \tau \displaystyle \int_{t-\eta}^{t} X(S) dS\Big)\\&+ r L_2 b_1 K_3 \Big(\displaystyle \int_{t-\delta}^{t}dS \displaystyle \int_{s}^{t} \lvert X(U-\tau) \rvert dU + \delta \displaystyle\int_{t-\tau}^{t} X(S)dS\Big)+ r b_1 L_2 K_4 \displaystyle \int_{t-\delta}^{t} dS\displaystyle \int_{s}^{t} \lvert Y(U) \rvert \\&+ r^{2} L_2^{2}\displaystyle \int_{t-\delta}^{t} dS \displaystyle \int_{s}^{t}\lvert Y(U) \rvert  + rL_2d_1\displaystyle\int_{t-\delta}^{t} dS \displaystyle \int_{s}^{t} \lvert Y(U) \rvert. 
\end{align*}
Therefore, we get;
\begin{align*}	
  &D^+ E_1 + D^+ E_2 \\&  \leq -\Big[(b K_1 - b_1 K_3 + d_1 + c M_1) - \eta (c b M_1 K_1 + c d M_1 + c^{2} M_1 M_2) - \tau(b_1 K_3^{2} b + b_1 d K_3\tau \\& + b_1 K_3 c M_2) - \delta (r L_2 b_1 K_3)\Big]\lvert X(t) \rvert - \Big[(b K_2 - b_1 K_4 + r L_1 - r L_2 +d_1) + (b_1 c M_1 K_2)\eta \\&- (b_1 b K_3K_4)\tau - (r L_2 b_1 K_4 + r^{2}L_2^{2} + 
 r L_2 d_1) \delta \Big] \lvert Y(t) \rvert - \Big[(c M_1\alpha)\eta - (b_1K_3\alpha)\tau\Big]\lvert Z(t) \rvert\\
& \leq - \Big[(A_1 - \eta B_1 - \tau C_1 - \delta D_1)\Big]\Big\lvert X(t) \Big\rvert - \Big[(A_2 - B_2 \eta - C_2 \eta - D_2 \delta) \Big]\Big\lvert Y(t) \Big\rvert - \Big[(B_3 \eta - C_3 \tau )\Big]\Big\lvert Z(t) \Big\rvert.
\end{align*}
 From our presumption, we have, $0<\eta<r$, $0<\tau<s$, $0<\delta<q$, and $D^+ E_1 + D^+ E_2<0$.\\[1mm]Hence, the equilibrium $(x^{*},y^{*},z^{*})$ is globally asymptotically stable for $0<\eta<r$, $0<\tau<s$ and $0<\delta<q$, where, $r=Min\Big\{\dfrac{A_1}{B_1}, \dfrac{A_2}{B_2}\Big\}$, $s=Min\Big\{\dfrac{A_1}{C_1}, \dfrac{A_2}{C_2},\dfrac{A_1 B_3}{C_3 B_1}, \dfrac{A_2 B_3}{C_3 B_2}\Big\}$ and $q=Min\Big\{\dfrac{A_1}{D_1}, \dfrac{A_2}{D_2}\Big\}$.
\end{proof}
\begin{remark} In contrast to local stability, which must abide by beginning states, global stability has been established to be independent of initial conditions. Thus, Theorem \ref{indep} gives parametric conditions for global stability, and Theorem \ref{dep} determines the length of the delay for which the solutions of the model reach equilibria globally. Basically, in this theorem, the length of the delay signifies the time period in which the contagious disease remains under control in a population. \\The delay-independent result derived for the model in \cite{rao2015dynamic} was for a particular case where $f(x,y) = xy$ (simple interaction function),$v(x) = x$ and $p(y) = y$ (linear functions), whereas the result Theorem \ref{indep} derived here is valid for the general case. Thus, it can be implemented for the model in \cite{rao2015dynamic}. Calculating the delay range for the global stability of the model is a new technique for considered infectious disease models. This methodology can also be applied to the model in \cite{rao2015dynamic} by taking $\eta=0$.
\end{remark}
Now, let's look at some numerical examples where we replicate the systems using various parametric values taken from COVID real-time data and observe their performance in order to validate our findings.
\section{Numerical examples and simulation}
\begin{figure}[!htb]
	\minipage{0.3\textwidth}
  \includegraphics[width = 5cm, height = 3.5cm]{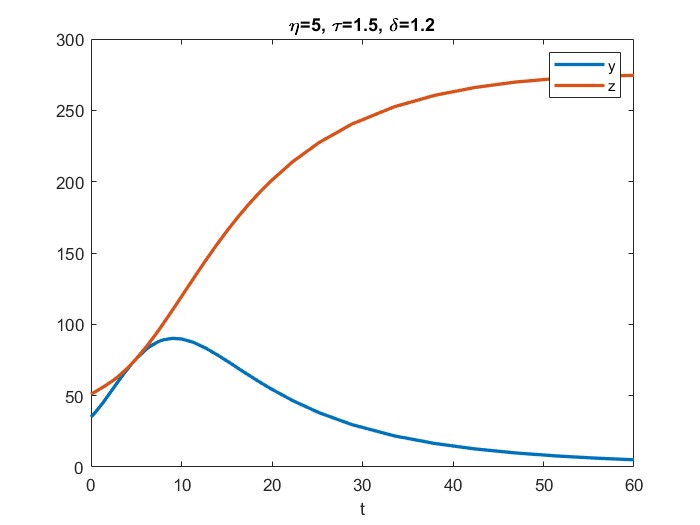}
  \caption{Solution curve at $\eta = 5, \tau = 1.5, \delta = 1.2$}
    \label{fig2}
\endminipage\hfill
	\minipage{0.3\textwidth}%
  \includegraphics[width = 5cm, height = 3.5cm]{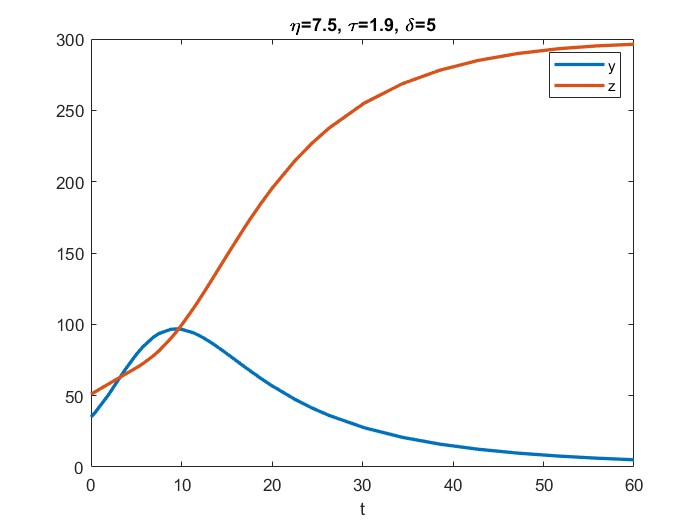}
  \caption{Solution curve at $\eta =  7.5, \tau = 1.9, \delta = 5$}
    \label{fig3}
\endminipage\hfill
\minipage{0.3\textwidth}
  \includegraphics[width = 5cm, height = 3.5cm]{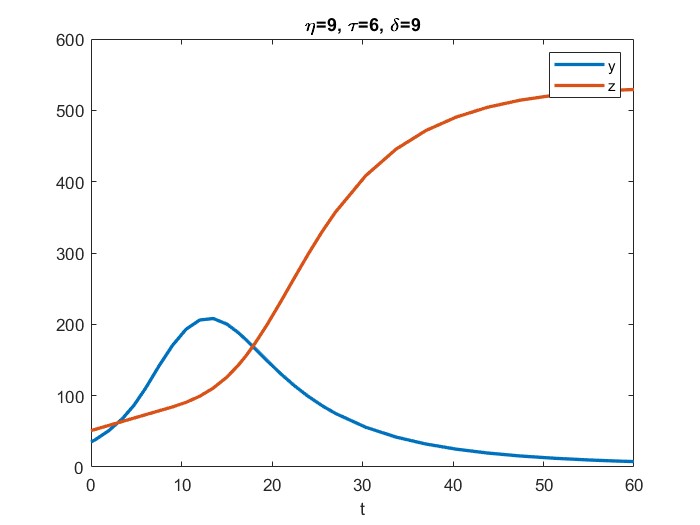}
  \caption{Solution curve at $\eta =  9, \tau = 6,\delta = 9$}
    \label{fig4}
\endminipage\hfill
\end{figure}

\begin{figure}[!htb]
	\minipage{0.3\textwidth}
  \includegraphics[width = 5cm, height = 3.5cm]{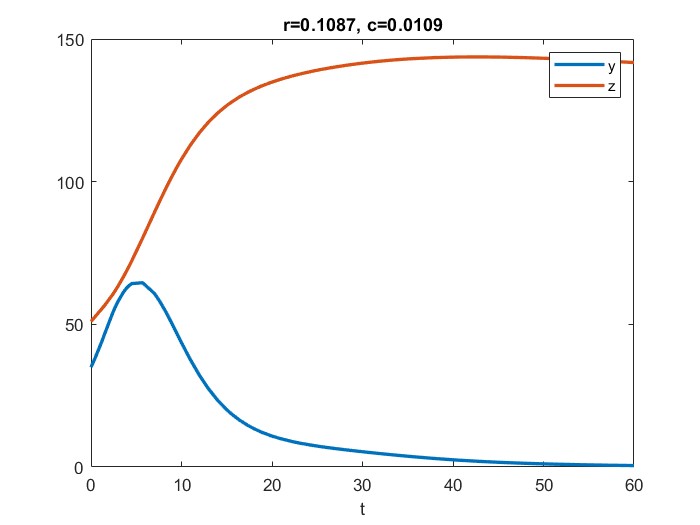}
  \caption{Solution curve at $r=0.1087, c=0.0109$}
    \label{fig5}
\endminipage\hfill
	\minipage{0.3\textwidth}%
  \includegraphics[width = 5cm, height = 3.5cm]{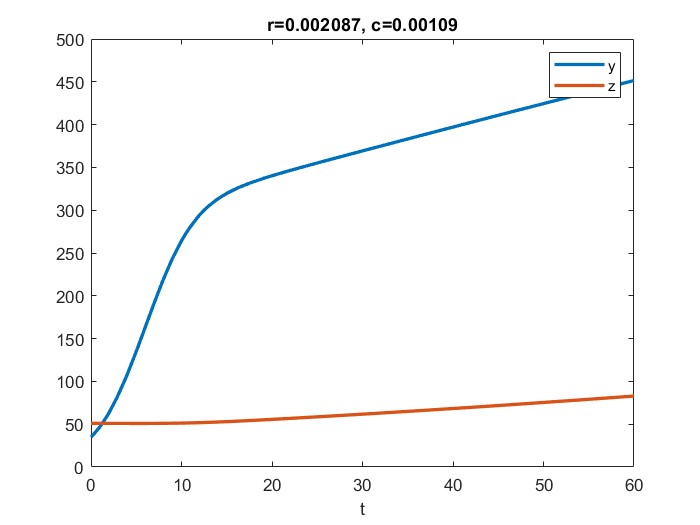}
  \caption{Solution curve at $r=0.002087, c=0.001090$}
    \label{fig6}
\endminipage\hfill
\minipage{0.3\textwidth}
  \includegraphics[width = 5cm, height = 3.5cm]{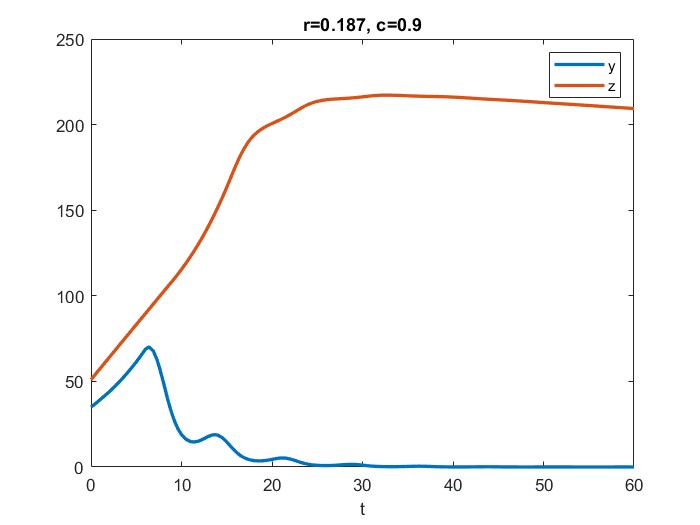}
  \caption{Solution curve at $r=0.1870, c=0.9$}
    \label{fig7}
\endminipage\hfill
\end{figure}

\begin{figure}[!htb]
\minipage{0.3\textwidth}
  \includegraphics[width = 5cm, height = 3.5cm]{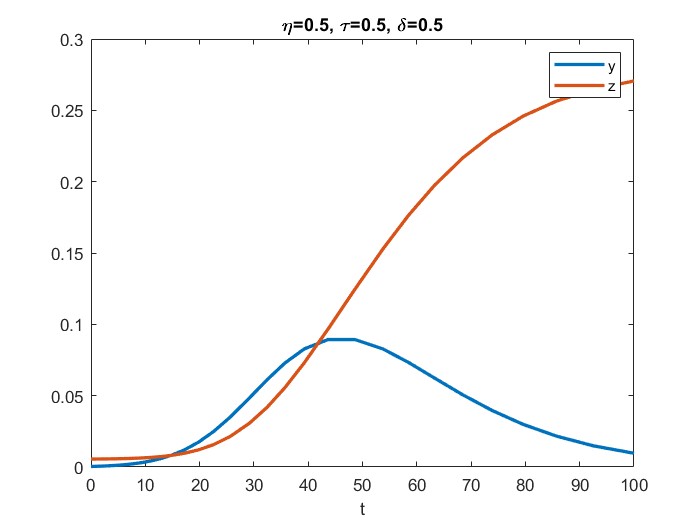}
  \caption{Solution curve at $\eta = 0.5, \tau = 0.5,  \delta = 0.5$}
  \label{fig14}
\endminipage\hfill
\minipage{0.3\textwidth}
  \includegraphics[width = 5cm, height = 3.5cm]{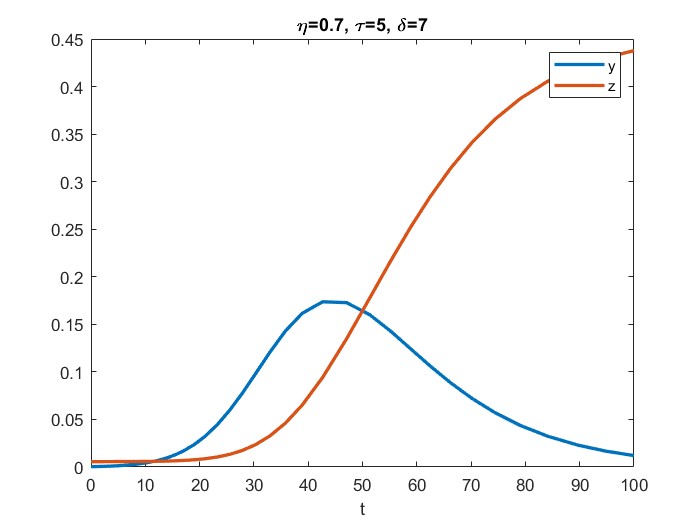}
  \caption{Solution curve at $\eta = 0.5, \tau = 5,  \delta = 7$}
  \label{fig15}
\endminipage\hfill
\minipage{0.3\textwidth}%
  \includegraphics[width = 5cm, height = 3.5cm]{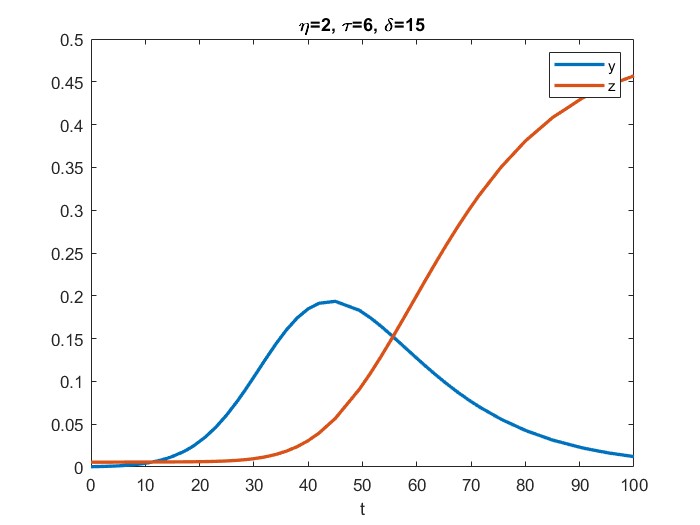}
  \caption{Solution curve at $\eta = 2, \tau = 6,  \delta = 15$}
  \label{fig16}
\endminipage\hfill
\end{figure}

\begin{figure}[!htb]
\minipage{0.3\textwidth}
  \includegraphics[width = 5cm, height = 3.5cm]{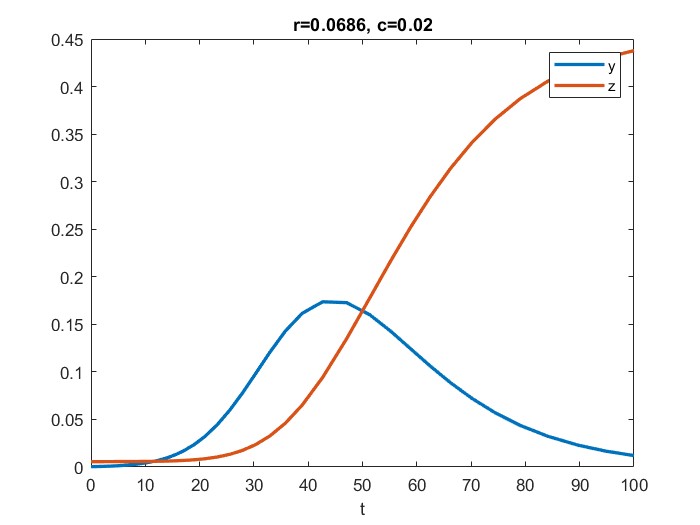}
  \caption{Solution curve at  $r = 0.0686, c=0.02$}
  \label{fig17}
\endminipage\hfill
\minipage{0.3\textwidth}%
  \includegraphics[width = 5cm, height = 3.5cm]{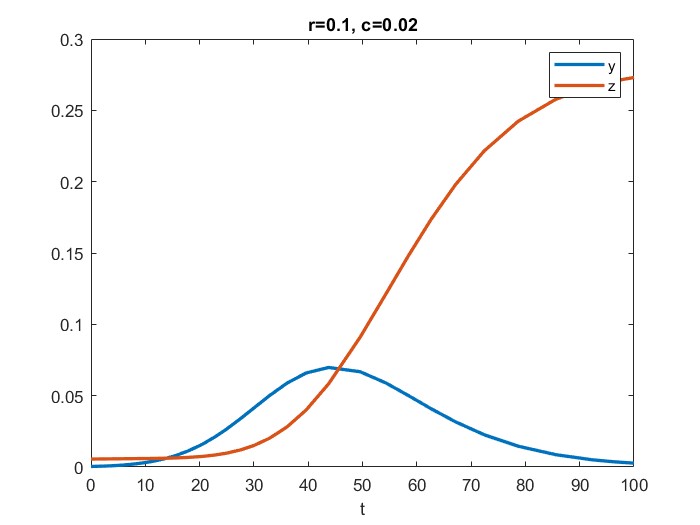}
  \caption{Solution curve at  $r = 0.1, c=0.02$}
  \label{fig18}
\endminipage\hfill
\minipage{0.3\textwidth}
  \includegraphics[width = 5cm, height = 3.5cm]{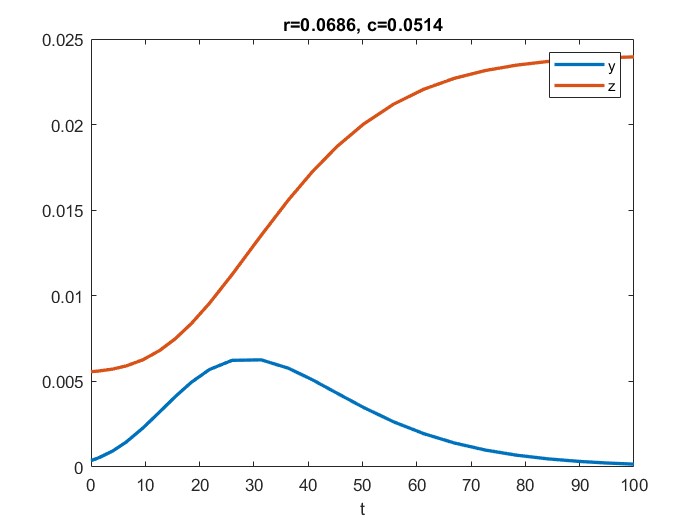}
  \caption{Solution curve at $r = 0.0686, c=0.0514$}
  \label{fig19}
\endminipage
\end{figure}

\begin{figure}[!htb]
\minipage{0.3\textwidth}
  \includegraphics[width = 5cm, height = 3.5cm]{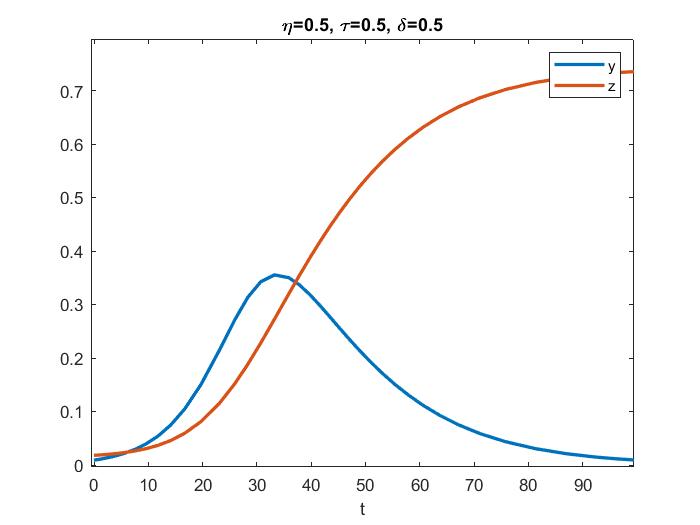}
  \caption{Solution curve at $\eta = 0.5, \tau = 0.5,  \delta = 0.5$}
  \label{fig8}
\endminipage\hfill
\minipage{0.3\textwidth}
  \includegraphics[width = 5cm, height = 3.5cm]{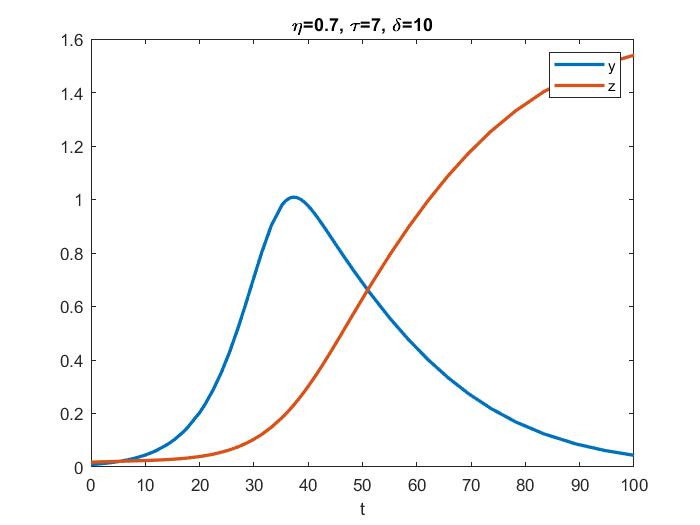}
  \caption{Solution curve at $\eta = 0.7, \tau = 7,  \delta = 10$}
  \label{fig9}
\endminipage\hfill
\minipage{0.3\textwidth}%
  \includegraphics[width = 5cm, height = 3.5cm]{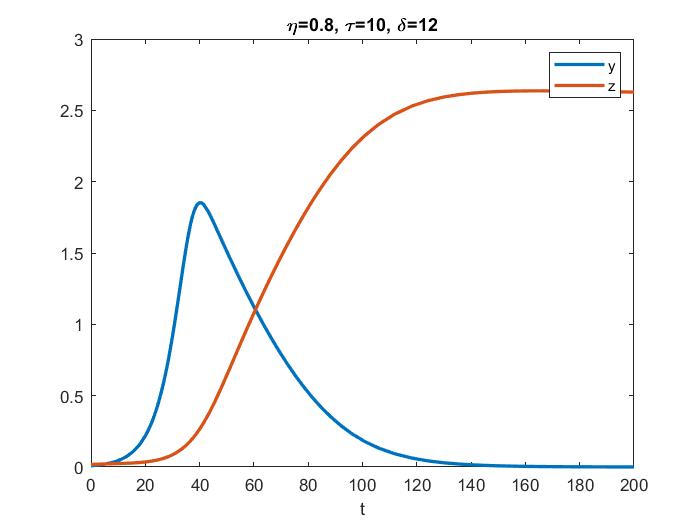}
  \caption{Solution curve at $\eta = 0.8, \tau = 10,  \delta = 12$}
  \label{fig10}
\endminipage
\end{figure}

\begin{figure}[!htb]
\minipage{0.3\textwidth}
  \includegraphics[width = 5cm, height = 3.5cm]{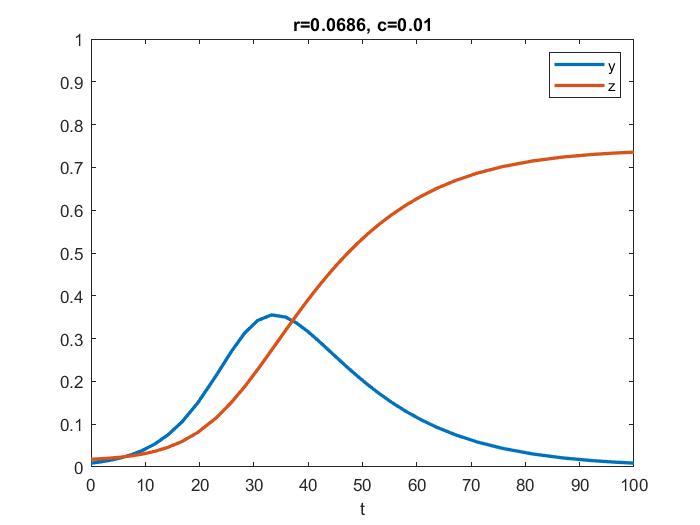}
  \caption{Solution curve at  $r = 0.0686, c=0.01$}
  \label{fig11}
\endminipage\hfill
\minipage{0.3\textwidth}%
  \includegraphics[width = 5cm, height = 3.5cm]{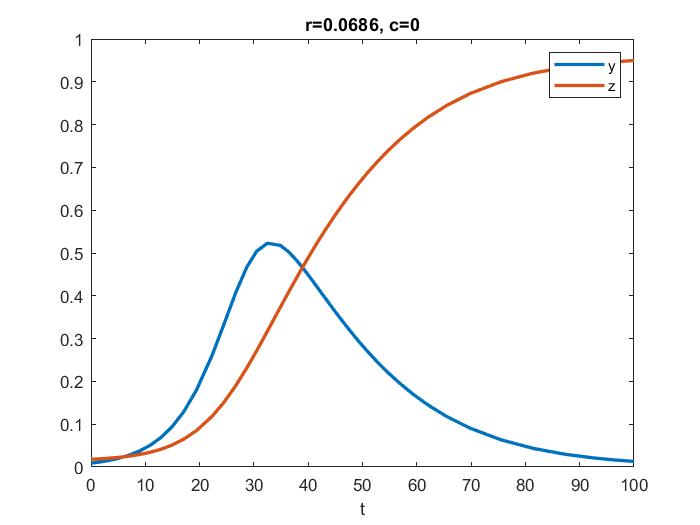}
  \caption{Solution curve at  $r = 0.0686, c=0$}
  \label{fig12}
\endminipage\hfill
\minipage{0.3\textwidth}
  \includegraphics[width = 5cm, height = 3.5cm]{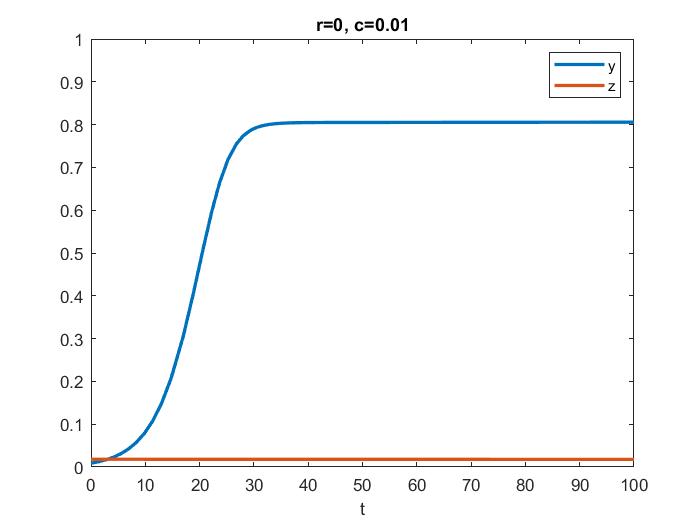}
  \caption{Solution curve at $r = 0, c=0.01$}
  \label{fig13}
\endminipage
\end{figure}

COVID-19 is a recent transmissible disease that has had a significant global impact. In this section, we adopt the parametric values from the real-time datasets of COVID-19 and examine the applicability of the model. To do so, we took three different sets of real-time data from COVID-19. The first dataset is the omicron variant of COVID-19 in the state of Tamilnadu, India (December 2021–March 2022); the second and third are the cumulative datasets of India and USA COVID-19 (November 2020), respectively. We shall set the parametric value from the above dataset and simulate the models by varying the delays to study their impact on disease prevalence and recovery. Additionally, we will make changes in the vaccination rates as well as the treatment rates to examine how the disease will eventually behave. These numerical examples have been simulated with the help of the MATLAB dde23 solver for systems of equations.\\[0.5mm]
Let us consider examples:
	\begin{example}
		\begin{align}
			x^{'} &=5-0.0012xy-0.065x-0.0109\Big(x(t-\eta)\Big)+0.0017z ,\nonumber \\
			y^{'} &=0.0012\Big(x(t-\tau)y\Big)-0.1087y-0.0006y,\nonumber\\
				z^{'} &=0.1087\Big(y(t-\delta)\Big)-0.0017z.\label{ex2}
		\end{align}
		\end{example}		
		In this example, we incorporate parametric values from the Tamilnadu COVID-19 dataset \cite{dickson2022sqirv}, therefore the parametric values are taken as $a = 5$, $b = 0.0012$, $c = 0.0109$, $d = 0.065$, $b_1 = 0.0012$, $r = 0.1087$, $d_1 = 0.0006$, $\alpha = 0.0017$ along with the infection function $f (x,y) = xy$, vaccination function $v(x) = x$ and treatment function $p(y) = y$. The basic reproduction number for this example is $R_0$ = $0.8344 < 1$. Therefore, by virtue of Theorem \ref{local2}, the system can be stable at equilibrium, provided other criteria of Theorem \ref{local2} must hold true. \\    Conditions of global stability in the context of the delay-dependent theorem \ref{dep} with delays range $\eta\in(0,13.6841)$,\hspace{0.1cm}$\tau\in (0,124.2974)$,\hspace{0.1cm}$\delta\in(0,0.050)$ are satisfied with initial conditions $(300, 35, 51)$. The significance of the below Figures \ref{fig2}, \ref{fig3}, and \ref{fig4} is how the infected population and recovered population change with respect to variation in delays. Figures \ref{fig5}, \ref{fig6}, and \ref{fig7} show how the change in the rate of vaccination and rate of treatment influences the infected and recovered populations.\\[2mm]
\textbf{Observations:}
From Figures \ref{fig2}, \ref{fig3}, \ref{fig4}, we have noted that as the duration of delays increases, infection increases, and as the duration of delays decreases, infection decreases. Figures \ref{fig5}, \ref{fig6}, \ref{fig7} show how the prevalence of disease declines as vaccination as well as recovery rates increase and rises when these rates are low.
\begin{example}
		\begin{align*}
			x^{'} &=0.00004893-0.4\Bigg(\dfrac{xy}{1+y}\Bigg)-0.00001992x-0.02\Big(x(t-\eta)\Big)+0.00017z , \\
			y^{'} &=0.32\Bigg(\dfrac{(x-\tau)y}{(x-\tau) +y}\Bigg)-0.0686\Bigg(\dfrac{y}{y + 1}\Bigg)-0.00002021y,\\
				z^{'} &=0.0686\Bigg(\dfrac{y(t-\delta)}{(y(t-\delta) + 1)}\Bigg)-0.00017z.
				\end{align*}
				\end{example}
				For this example, we took India's COVID-19 real-time parametric values, referred to as \cite{wintachai2021stability} where, $a = 0.00004893$, $b = 0.4$, $c = 0.02$, $d = 0.00001992$, $b_1 = 0.32$, $r = 0.0686$, $d_1 = 0.00002021$, $\alpha = 0.00017$.  Here, we changed the infection function as $f(x,y) = \dfrac{xy}{(1 + y)}$, treatment function $p(y) = \dfrac{y}{y+1}$ and vaccination function $v(x) = x$,  as semi-linear and linear, respectively. We try to put various inputs and tracked the nature of the solution curve, which is shown in Figures \ref{fig14}, \ref{fig15}, \ref{fig16}, \ref{fig17}, \ref{fig18}, \ref{fig19}.        The basic reproduction number value for this specific example is $R_0$ = $11.4545 > 1$, which reveals that the system is locally unstable. The delay range for the condition of global stability (delay-dependent) is given by $\eta\in(0,11.9066)$,\hspace{0.1cm}$\tau\in (0,0.6252)$,\hspace{0.1cm}$\delta\in(0,3.0016)$ are satisfied with the given initial conditions $(0.994, .0003813, .005569)$.\\[2mm]
\textbf{Observations:}
As in the above example, Figures \ref{fig14}, \ref{fig15}, \ref{fig16} depict the behavior of disease as the delays vary. After a comparative study of Figures \ref{fig17}, \ref{fig18}, \ref{fig19}, we noticed that if only the vaccination rate or treatment rate increases, even then there is a chance to control the effects of the disease in society. It is also observed that the rate of recovery is higher when there is an increase in the vaccination rate compared to an increase in the treatment rate.
	\begin{example}
		\begin{align*}
			x^{'} &=0.000031785-0.5\Bigg(\dfrac{xy}{1 + x}\Bigg)-0.00002377x-0.01\Big(x(t-\eta)\Big)+0.00017z , \\
			y^{'} &=0.462\Bigg(\dfrac{x(t - \tau)y}{(1 + x(t - \tau)}\Bigg)-0.0686\Bigg(\dfrac{y}{y + 1}\Bigg)-0.00002585y,\\
				z^{'} &=0.0686\Bigg(\dfrac{y(t-\delta)}{(y(t-\delta) + 1)}\Bigg)-0.00017z.
		\end{align*}
		\end{example}	
		In this particular example, we initiate with US COVID-19 real-time parametric values taken from source \cite{wintachai2021stability} as\hspace{0.1cm} $a = 0.000031785$, $b = 0.5$, $c = 0.01$, $d = 0.00002377$, $b_1 = 0.462$, $r = 0.0686$, $d_1 = 0.00002585$, $\alpha = 0.00017$. We set the infection function as  $f (x,y) = \dfrac{xy}{(1 + x)}$, treatment function $p(y) = \dfrac{y}{y+1}$, and vaccination function $v(x) = x$. In this particular instance, the basic reproduction number value is  $R_0$ = $3.8517 > 1$, which infers that the system is locally unstable. Global stability conditions in the frame of delay-dependent propositions involving delays ranging from $\eta\in(0,0.8.2307)$,\hspace{0.1cm}$\tau\in (0,0.1646)$,\hspace{0.1cm}$\delta\in(0,1.0446)$ are satisfied with the given initial conditions $(0.97286, 0.00905, 0.01809)$. We observed the importance of Figures \ref{fig8}, \ref{fig9}, \ref{fig10},\ref{fig11}, \ref{fig12}, \ref{fig13} in view of the epidemic. \\[2mm]
\textbf{Observations:}
According to our observation from Figures \ref{fig8}, \ref{fig9}, \ref{fig10}, as the duration of delays rises, infection rises as well, and vice versa as the duration of delays falls. On comparing Figure \ref{fig11}, \ref{fig12}, and \ref{fig13}, it is observed in Figure \ref{fig12} that, in the absence of vaccination, there will be a spike in the number of infected people, but because of the availability of treatment, disease prevalence will be under control. Whereas in Figure \ref{fig13}, it is noted that even when there is the availability of vaccination, if there is no treatment, disease prevails continuously in the population.    \begin{remark}In the above-mentioned examples, we have verified our model with real-time data. In the absence of treatment, vaccination efforts are inadequate to manage the conditions. Treatment is advised in both the severe and long-term stages of the disease. Our examples lead us to conclude that an environment free of disease arises when treatment and vaccination rates are high. On the other hand, if both the treatment and vaccination rates are low, the disease will predominate.
	   
	 \end{remark}
	\section{Conclusions and future work}
	 In this article, we explore an extended SIR model that accounts for the time delays in the functions of infection, vaccination, and treatment. We employed general non-linear infection functions to analyze the dynamics of disease transmission in the human population. In addition to other qualitative aspects like non-negativity and boundedness, the existence-uniqueness of solutions and equilibrium, the basic reproduction number, has been discussed because of its biological significance, in terms of understanding the dynamics of a disease in a population. The basic reproduction number plays a significant role in the model's local stability.
	The model that is being provided here is rather unique in the sense of global stability in terms of delay dependence, where the range of the delays has been established along with a few parametric constraints under the implications of an appropriate Lyapunov function. Another factor that emphasizes the significance of this study is how treatment and vaccination work together. It has been found that there is a permissible vaccination threshold rate at which the pathogen is completely eradicated. Numerical examples and simulations based on various real-time datasets indicate how the behavior of the solutions changes when vaccination, infection, and treatment functions act together under the influence of delays.\\[0.6mm]Even though the fact that there are many SIR models given in the literature. However, from the perspective of future scope, there is always a chance for modification in those models. Many researchers introduced model stability modulated by vaccination and treatment function in the absence of delays \cite{saha2022dynamical,zaman2017optimal} but in this paper, we have covered these aspects by introducing the delays in these functions. In our model assumptions, we focused on the susceptible population's growth rate, which is assumed to be a fixed constant, and this presumption is valid when a significant amount of people contributes to the population. In the case of a disease like AIDS, each one of these compartments could make a substantial contribution to the system, so we need to explore the dynamic behavior of the disease with an adjustable population. 
	In addition to this, discretization techniques and the stability behavior of dynamical systems under the influence of non-linear functions of continuous models are also still unpredictable \cite{rao2015dynamic}. Real-time data fitting is one of the major challenges in recent research. The utility of mathematical models is only acceptable if they pass the test on real-world data, even though they can describe a wide range of plausible and realistic phenomena. Data fitting will make the model more realistic, along with parametric values. Further, we can also extend this model into more compartments for better analysis of an epidemic. Therefore, it would be very interesting to test the outcomes of this study using more appropriate real-time data so that any enhancements or modifications might be adequately implemented.\\\\
\textbf{Declaration of competing interest.}\\
The authors declare that they have no known competing financial interests or personal relationships that could
have appeared to influence the work reported in this paper.\\\\
\textbf{Author contributions.}\\ The authors contributed equally to the completion of this research work.

	\bibliographystyle{plain}	
	\bibliography{Citations} 
\end{document}